\documentclass[12pt]{amsart}
\usepackage{amssymb,amsmath}

\makeatletter
\def\@cite#1#2{{\m@th\upshape\bfseries%
[{#1\if@tempswa{\m@th\upshape\mdseries, #2}\fi}]}}
\makeatother
\let\secsymb=\S
\theoremstyle{plain}
\newtheorem{thm}{Theorem}[section]
\newtheorem{cor}[thm]{Corollary}

\newtheorem{lem}[thm]{Lemma}

\theoremstyle{definition}

\newtheorem{eg}[thm]{Example}

\newcommand{\bB}{{\mathbb{B}}}
\newcommand{\bC}{{\mathbb{C}}}
\newcommand{\bD}{{\mathbb{D}}}

\newcommand{\bF}{{\mathbb{F}}}

\newcommand{\bN}{{\mathbb{N}}}

\newcommand{\bS}{{\mathbb{S}}}

 \newcommand{\A}{{\mathcal{A}}}
 \newcommand{\B}{{\mathcal{B}}}

 \newcommand{\E}{{\mathcal{E}}}
 \newcommand{\F}{{\mathcal{F}}}
 
\renewcommand{\H}{{\mathcal{H}}}

 \newcommand{\K}{{\mathcal{K}}}

 \newcommand{\M}{{\mathcal{M}}}
 \newcommand{\N}{{\mathcal{N}}}

\renewcommand{\P}{{\mathcal{P}}}

\renewcommand{\S}{{\mathcal{S}}}
 
 \newcommand{\U}{{\mathcal{U}}}

\newcommand{\eps}{\varepsilon}
\renewcommand{\phi}{\varphi}

\newcommand{\upchi}{{\raise.35ex\hbox{$\chi$}}}

\newcommand{\fL}{{\mathfrak{L}}}
\newcommand{\fM}{{\mathfrak{M}}}

\newcommand{\fR}{{\mathfrak{R}}}

\newcommand{\rC}{\mathrm{C}}

\newcommand{\ca}{\mathrm{C}^*}

\newcommand{\Fn}{\bF_n^+}
\newcommand{\Focksp}{\ell^2(\Fn)}
\newcommand{\lip}{\langle}
\newcommand{\rip}{\rangle}
\newcommand{\ip}[1]{\left\lip #1 \right\rip}

\newcommand{\mt}{\varnothing}

\newcommand{\ol}{\overline}

\newenvironment{sbmatrix}{\left[\begin{smallmatrix}}{\end{smallmatrix}\right]}

\newcommand{\upminus}{\raise.9ex\hbox{-\!}}

\newcommand{\wot}{\textsc{wot}}

\newcommand{\AD}{\mathrm{A}(\mathbb{D})}

\newcommand{\CT}{{\mathrm{C}(\mathbb{T})}}

\newcommand{\Fock}{{\F_n}}
\newcommand{\sFock}{{\H_n^2}}
\newcommand{\aFock}{{\A_n}}

\newcommand{\FORAL}{\text{ for all }}

\newcommand{\qand}{\quad\text{and}\quad}

\newcommand{\qfor}{\quad\text{for}\quad}

\newcommand{\diam}{\operatorname{diam}}

\newcommand{\Mult}{\operatorname{Mult}}

\begin{document}

\title[Commutant Lifting for Commuting Row Contractions]%
{Commutant Lifting\\ for Commuting Row Contractions}
%
\author[K.R.Davidson]{Kenneth R. Davidson}
\thanks{First author partially supported by an NSERC grant}
\address{Pure Math.\ Dept.\\U. Waterloo\\Waterloo, ON\;
N2L--3G1\\CANADA}
\email{krdavids@uwaterloo.ca}
\author{Trieu Le}
\address{Pure Math.\ Dept.\\U. Waterloo\\Waterloo, ON\;
N2L--3G1\\CANADA}
\email{t29le@uwaterloo.ca}
\subjclass[2000]{47A20}
\date{}

\begin{abstract}
If $T= \big[ T_1\ \dots\ T_n\big]$ is a row contraction with commuting entries,
and the Arveson dilation is $\tilde T= \big[ \tilde T_1\ \dots\ \tilde T_n\big]$, 
then any operator $X$ commuting with each $T_i$
dilates to an operator $Y$ of the same norm which commutes with each $\tilde T_i$.
\end{abstract}
\maketitle

\section{Introduction}\label{S:intro}

The commutant lifting theorem of Sz.Nagy and Foia\c{s} \cite{SF_CLT,SF}
is a central result in the dilation theory of a single contraction.  
It states that if $T \in \B(\H)$ is a contraction with isometric dilation 
$V$ acting on $\K \supset \H$, and $TX=XT$,
then there is an operator $Y$ with $\|Y\|=\|X\|$, $VY=YV$ and $P_\H Y = XP_\H$.
This result is equivalent to Ando's Theorem that two commuting contractions have
a joint (power) dilation to commuting isometries.
However the Ando dilation is not unique, and Varopoulos \cite{Var} showed that
it fails for a triple of commuting contractions.
In particular, the commutant lifting result does not generalize if one replaces 
$T$ by a commuting pair $T_1$ and $T_2$ \cite{Mull}.
See Paulsen's book \cite{Pau} for a nice treatment of these issues.

There is a multivariable context where the commutant lifting theorem does hold.
This is the case of a row contraction $T= \big[ T_1\ \dots\ T_n \big]$, i.e.\ $\|T\|\le1$
considered as an operator in $\B(\H^{(n)},\H)$.  
The operators $T_i$ are not assumed to commute.
The Frazho--Bunce dilation theorem \cite{Fr,Bun} shows that there is a
minimal dilation to a row isometry.
Popescu \cite{Pop1} establishes uniqueness of this dilation.
Frazho \cite{Fr2} and Popescu \cite{Pop1} establish the corresponding commutant lifting theorem.
It works because the proofs of Bunce and Popescu follow the original 
Sz.Nagy--Foia\c{s} proofs almost verbatim.

There is a related multivariable context of commuting row contractions.
This area has received a lot of interest recently beginning with several
papers of Arveson \cite{Arv3,Arv4,Arv5}.
The analogue of the von Neumann inequality in this context was established 
by Drury \cite{Drury}.  The extension to an appropriate dilation theorem is due
to Muller--Vasilescu \cite{MV} and Arveson \cite{Arv3}.
We follow Arveson's approach.  The model row contraction is the $n$-tuple
$M= \big[ M_{z_1}\ \dots\ M_{z_n} \big]$ of multipliers on symmetric Fock space $\sFock$,
sometimes called the Drury--Arveson space.
Every commuting row contraction $T$ dilates to a canonical operator 
$\tilde{T} \simeq M^{(\alpha)} \oplus U$, where $U$ is a spherical unitary, i.e.\ 
$U = \big[ U_1\ \dots\ U_n \big]$ where $U_i$ are commuting normal operators
satisfying $\sum U_iU_i^*=I$. (See the next section for details.)

A commutant lifting theorem in a related context was established by Ball, Trent and
Vinnikov \cite{BTV}.  They only consider $n$-tuples with a pure dilation; that is,
$T$ dilates to some $M^{(\alpha)}$ with no spherical unitary part.
They work more generally in the framework of complete Nevanlinna--Pick kernel
Hilbert spaces.  
The connection is made because of the characterization of 
complete Nevanlinna--Pick kernels by McCullough \cite{McC1,McC2}, 
Quiggin \cite{Q} and Agler--McCarthy \cite{AM}.
In the last section of this paper, we will discuss their result
and show how it follows from ours. 
Popescu \cite{Pop_reln} has a generalization of the Ball--Trent--Vinnikov result
for pure row contractions satisfying relations other than commutativity.

In this paper, we establish the commutant lifting theorem for commuting row contractions.

\begin{thm} \label{CLT}
Suppose that $T=\big[ T_1\ \dots\ T_n \big]$ is a commuting  row contraction 
on a Hilbert space $\H$, and that $X$ is an operator on $\H$ which commutes 
with $T_1,\dots, T_n$. 
Let $\tilde T = \big[ \tilde T_1\ \dots\ \tilde T_n \big]$ be the Arveson dilation 
of $T$ on $\K$. 
Then there is an operator $Z$ on $\K$ that 
commutes with each $\tilde{T}_i$ for $1\le i\le n$, 
which dilates $X$ in the sense $P_\H Z = X P_\H$
and satisfies $\|Z\|=\|X\|$.
\end{thm}

The minimal isometric dilation of a contraction is actually a co-extension.
That is, given a contraction $T$, the minimal isometric dilation is an isometry
$V$ on a Hilbert space $\K \supset \H$ such that $\H$ is co-invariant for $V$ and
$P_\H V|_\H = T$.  It follows that $P_\H V^k|_\H = T^k$ for all $k \ge 0$.
The Arveson dilation is also a co-extension.

The minimal unitary dilation $U$ of a contraction $T$ is of a more general type.
It is still true that $P_\H U^k|_\H = T^k$ for all $k \ge 0$.
This implies \cite{Sar} that $\H$ is semi-invariant, and we can decompose $\K$ into a direct sum 
$\K = \H_-\oplus\H\oplus \H_+$ so that $U$ has the form
$U = \begin{sbmatrix} * & 0 & 0 \\ * & T & 0 \\ * & * & * \end{sbmatrix}$.
The restriction of $U$ to $\H\oplus \H_+$ is the minimal isometric dilation $V$.
The commutant lifting theorem of Sz.Nagy and Foia\c{s} is still valid for
this unitary dilation, or indeed, for any unitary dilation of $T$.

In the case of a commuting row contraction, the Arveson dilation is unique,
but a maximal dilation of a row contraction of the type considered in the 
previous paragraph is not.  
However Arveson \cite{Arv3} showed that the C*-envelope of the 
universal commuting row contraction
is just $\ca(M) = \ca(\{M_{z_1},\dots,M_{z_n}\})$.
The  Arveson dilation is already maximal in the sense that any further dilation
(either an extension or co-extension) to a commuting row contraction will 
just be the addition of a direct summand of another commuting row contraction.
Therefore the Arveson dilation is canonical and plays a special role,
and the more general dilations are not really needed.

Nevertheless, one can ask whether the commutant lifting theorem is valid if one allows
these more general dilations.  It should not be surprising that the lack of
uniqueness makes such a result impossible, as an easy example shows.

We add that there is a context in which one can deal with lack of uniqueness
of unitary dilations and still retain commutant lifting.
The tree algebras of Davidson--Paulsen--Power \cite{DPP} have such a result.
Given a completely contractive representation commuting with a contraction $X$,
one can find \textit{some} unitary dilation of the algebra commuting with a
contractive dilation of $X$.  If one restricts one's attention to the smaller class  
that Muhly and Solel \cite{MS} call tree algebras, the minimal unitary dilation
is unique and one obtains the more traditional form of commutant lifting.

\section{Background}\label{S:back}

We start with the canonical model for a row contraction.
Let $\Fn$ denote the free semigroup on $n$ letters, and form 
Fock space $\Fock = \Focksp$ with basis $\{ \xi_w : w \in \Fn \}$. 
Define the left regular representation of $\Fn$ on $\Fock$ by
$L_v \xi_w = \xi_{vw}$ for $v,w \in \Fn$.
The algebra $\fL_n$ is the \wot-closed span of $\{L_v : v \in \Fn \}$.
The commutant of $\fL_n$ is $\fR_n$, the algebra generated by the
right regular representation \cite[Theorem~1.2]{DP1}.

In particular, $L = \big[ L_1\ \dots \ L_n \big]$ is a row isometry
with range $(\bC \xi _\mt)^\perp$.
This row isometry $L$ is the canonical model for a row contraction
because of Popescu's von Neumann inequality \cite{Pop2}
which states that for any $A = \big[ A_1\ \dots \ A_n \big]$
which is a row contraction, i.e.\ $\|A\| \le 1$, and any polynomial
$p$ in $n$ non-commuting variables, one has
\[ \|p(A_1,\dots,A_n) \| \le \| p(L_1,\dots,L_n) \| .\]

Now we consider the commuting case.
Let $\bN_0$ be the set of non-negative integers.  
For $\lambda\in\bC^n$ and each
$n$-tuple $k=(k_1,\dots,k_n)\in\bN_0^n$,  let us write
$\lambda^k:=\lambda_1^{k_1}\dots\lambda_n^{k_n}$.  
For $k\in\bN_0^n$, let $\P_k:=\{w\in\F_n:
w(\lambda)=\lambda^k \FORAL \lambda\in\bC^n\}$.  
Define vectors in $\Fock$ by
\[
  \zeta^k := \frac{1}{|\P_k|} \sum_{w\in \P_k} \xi_w .
\]
Note that $ |\P_k| =  \binom{|k|!}{k_1!\,k_2!\cdots k_n!}$
and $\|\zeta^k\|=|\P_k|^{-1/2}$.
This set of symmetric words forms an orthogonal basis
for \textit{symmetric Fock space} $\sFock$.
We consider $\sFock$ as a space of analytic functions
on the unit ball $\bB_n$ of $\bC^n$ by identifying
an element $f = \sum_{k\in\bN^n} a_k\zeta^k$ with the function
\[ f(\lambda) = \sum_{k\in\bN^n} a_k \lambda^k .\]
This series converges uniformly for $\|\lambda\| \le r$ for any $r<1$.

For each $\lambda \in \bB_n$, let 
\[
 k_\lambda =  \sum_{w\in\F_n} \ol{w(\lambda)} \xi_w  =
 \sum_{k\in\bN_0^n} \ol{\lambda}^k |\P_k| \zeta^k .
\]
Then 
\[ \ip{\zeta^k, k_\lambda} = \lambda^k |\P_k| \|\zeta^k\|^2 = \lambda^k .\]
Therefore $\ip{f,k_\lambda} = f(\lambda)$ for all $\lambda\in\bB_n$.
It is evident that each $k_\lambda$ belongs to $\sFock$; and they
span the whole space because any element $f\in\sFock$ orthogonal to 
all $k_\lambda$ satisfies $f(\lambda) = 0$ for $\lambda\in\bB_n$, and hence $f=0$.
Therefore $\sFock$ becomes a reproducing kernel Hilbert space.

A \textit{multiplier} on $\sFock$ is a function $h$ on $\bB_n$ so that 
$M_h f = hf$ determines a well defined map of $\sFock$ into itself.
A standard argument shows that such maps are continuous.
It is easy to see that $h$ must be an analytic function, and the set
of all multipliers forms a \wot-closed algebra of operators.
A routine calculation shows that $M_h^* k_\lambda = \ol{h(\lambda)} k_\lambda$
for $\lambda\in\bB_n$; and conversely, any bounded operator $T$ such
that $T^* k_\lambda = \ol{h(\lambda)} k_\lambda$ determines a multiplier $h$.

Now, \cite[Theorem~2.6]{DP1} shows that the vectors $k_\lambda$ 
are precisely the set of eigenvectors for $\fL_n^*$, and
$L_i^* k_\lambda = \ol{\lambda_i} k_\lambda$.
In particular,  
\[ P_\sFock L_i|_\sFock = M_{z_i} \qfor 1 \le i \le n \]
are the multipliers  by the coordinate functions.
Indeed, $\sFock$ is co-invariant for $\fL_n$ and the compression
to $\sFock$ yields a complete quotient of $\fL_n$ by its commutator ideal
onto the space of all multipliers of $\sFock$ \cite[Corollary~2.3 and \secsymb~4]{DP3}.
Here $A \in \fL_n$ is sent to $P_\sFock A|_\sFock = M_{\hat A}$ where
$\hat A(\lambda) := \ip{A k_\lambda,k_\lambda}/\|k_\lambda\|^2$.

The row operator $M = \big[ M_{z_1}\ \dots \ M_{z_n} \big]$ is the 
canonical model for a commuting row contraction.
This is because of Drury's von Neumann inequality \cite{Drury}
which states that for any $T = \big[ T_1\ \dots \ T_n \big]$
which is a commuting row contraction, i.e.\ $\|T\| \le 1$ and 
$T_iT_j=T_jT_i$ for $1 \le i,j \le n$, and any polynomial
$p$ in $n$ commuting variables, one has
\[ \|p(T_1,\dots,T_n) \| \le \| p(M_{z_1},\dots,M_{z_n}) \| .\]

If $A$ is an operator on a Hilbert space $\H$, we
call an operator $B$ on a Hilbert space $\K \supset \H$ an 
\textit{extension} of $A$ if $\H$ is invariant for
$B$ and $B|_\H =A$, i.e., $B =
\begin{sbmatrix} A & * \\ 0 & * \end{sbmatrix}$
with respect to the decomposition $\K = \H \oplus \H^\perp$.
Likewise, $B$ is a \textit{co-extension} of $A$ if $\H$ is
co-invariant for $B$ and $P_\H B = A P_\H$, i.e.\ $B = 
\begin{sbmatrix} A & 0 \\ * & * \end{sbmatrix}$.
Finally, we say that $B$ is a \textit{dilation} of $A$ if
$\H$ is semi-invariant and $P_\H B|_\H = A$, i.e.\ $B = 
\begin{sbmatrix}*&*&*\\ 0& A & * \\ 0 & 0 & * \end{sbmatrix}$
with respect to a decomposition $\K = \K_+ \oplus \H \oplus \K_-$.
By a result of Sarason \cite{Sar}, this latter spatial condition
is equivalent to the algebraic condition $P_\H B^k|_\H = A^k$
for $k \ge 0$.  This is sometimes called a power dilation.

The Sz.~Nagy dilation theory shows that every contraction $A$ 
has a unique minimal isometric co-extension, and a unique minimal
unitary dilation.
There is an analogue of the Sz.~Nagy isometric dilation in both 
of the multivariable contexts given above.
In the case of a row contraction $A = \big[ A_1\ \dots \ A_n \big]$, 
the Frazho--Bunce dilation \cite{Fr,Bun} states that there is a
unique minimal isometric co-extension of $A$.  That is, there is
a unique minimal row isometry $W =  \big[ W_1\ \dots \ W_n \big]$
on a Hilbert space $\K \supset \H$ so that each $W_i$ is a co-extension
of $A_i$.  Bunce's proof is to apply Sz.~Nagy's isometric dilation to
the row operator $A$ and observe that this suffices.

A row isometry $W$ on $\H$ must consist of $n$ isometries with
pairwise orthogonal ranges, and thus satisfies the Cuntz--Toeplitz relation 
$\sum_{i=1}^n W_iW_i^* \le I$. 
It is said to be Cuntz type if $\sum_{i=1}^n W_iW_i^* = I$.
It can be canonically decomposed using the Wold decomposition
\cite{Pop1} into $W \simeq L^{(\alpha)} \oplus S$ on 
$\H = \Fock^{(\alpha)} \oplus \K$, where $S$ is a Cuntz row isometry,
simply by setting $\M$ to be the range of $I-\sum_{i=1}^n W_iW_i^*$,
identifying the invariant subspace it generates with $\M \otimes \Fock$,
and letting $\K$ be the orthogonal complement.
Popescu's von Neumann inequality is a consequence of
the Frazho--Bunce dilation and the fact that there is a canonical
$*$-homomorphism from the Cuntz--Toeplitz C*-algebra 
$\E_n = \ca(L_1,\dots,L_n)$ onto the
$\ca(W_1,\dots,W_n)$ for any row isometry $W$.

In the case of a commuting row contraction, Drury did not extend his
von Neumann inequality to a dilation theorem.  This was done by
Muller and Vasilescu \cite{MV} in a very general but combinatorial way,
and by Arveson \cite{Arv3} in the context of symmetric Fock space.
He first notes that $\ca(M_{z_1}, \dots, M_{z_n})$ contains the compact
operators, and the quotient is $\rC(\bS_n)$, the space of continuous functions
on the unit sphere $\bS_n$ in $\bC^n$.
A \textit{spherical unitary} is a row operator $U= \big[ U_1\ \dots \ U_n \big]$
where the $U_i$ are commuting \text{normal} operators such that
$\sum_{i=1}^n U_iU_i^* = I$, i.e.\ the $U_i$ are the images of $z_i$ for $1 \le i \le n$
under a $*$-representation of $\rC(\bS_n)$.
Arveson's dilation theorem states that any commuting row isometry
$T = \big[ T_1\ \dots \ T_n \big]$ has a unique minimal co-extension
to an operator of the form $M^{(\alpha)} \oplus U$, 
where $M = \big[ M_{z_1}\ \dots \ M_{z_n} \big]$ is the multiplier on 
symmetric Fock space and $U$ is a spherical unitary.

Next, we consider the commutant lifting theorem.  
Consider the case of a single contraction $A$ on $\H$ with 
minimal isometric dilation $V$ on $\K$ and 
minimal unitary dilation $U$ on $\K'$.
Any contraction $X$ commuting with $A$ has a co-extension to 
a contraction $Y$ on $\K$ which commutes with $V$, and 
this has an extension to a contraction $Z$ on $\K'$ commuting with $U$.

Popescu \cite{Pop1} observes that exactly the same proof yields
a commutant lifting theorem for a row contraction $A = \big[ A_1\ \dots \ A_n \big]$
with minimal row isometric dilation $W$.
If there is a contraction $X$ which commutes with each $A_i$, then
there is a contraction $Y$ which is a co-extension of $X$  that 
commutes with $W_1,\dots,W_n$.
Frazho \cite{Fr2} did it for the case $n=2$.

In \cite{DP3}, the first author and Pitts gave a new proof of Drury's
von Neumann inequality by observing that when $T$ is a commuting
row contraction of norm less than 1 (so that the isometric dilation
has the form $L^{(\alpha)}$), the row isometric dilation on $\Fock^{(\alpha)}$
identifies the original space $\H$ with a subspace of $\sFock^{(\alpha)}$.

We have a similar approach to the commutant lifting theorem for
a commuting row contraction.  The plan is to apply Popescu's
commutant lifting theorem and restrict down to the space on which
Arveson's dilation lives.

\section{Commutant Lifting}\label{S:CLT}

To set the stage, let $T= \big[ T_1\ \dots\ T_n\big]$ be a commuting row contraction
on a Hilbert space $\H$, and let $X$ be an operator of norm 1 which commutes with 
each $T_i$, $1 \le i \le n$.
Let the Arveson dilation of $T$ be $\tilde T= \big[ \tilde T_1\ \dots\ \tilde T_n\big]$.
Decompose $\tilde T \simeq M^{(\alpha)} \oplus U$, where $U$ is a spherical unitary,
on $\K \simeq \sFock^{(\alpha)} \oplus \N$.
By the Frazho--Bunce dilation theorem, we can dilate $\tilde T$ to a row isometry 
$\hat T$ on $\K' \simeq \Fock^{(\alpha)} \oplus \N'$ which is unitarily equivalent to 
$L^{(\alpha)} \oplus S$, where $S$ is the isometric dilation of $U$ 
and is of Cuntz type because $\sum_{i=1}^n U_iU_i^*=I$.

The starting point of our proof is to invoke Popescu's commutant lifting theorem
to obtain an operator $Y$ on $\K$ commuting with $\hat T_i$ for $1 \le i \le n$
such that $\|Y\|=1$ and $P_\H Y = X P_\H$, so $Y$ is a co-extension of $X$.
Define $Z = P_\K Y|_\K$.  We will establish that $Z$ is the desired lifting.


The first lemma deals with the structure of the commutant of a row isometry
which dilates $M^{(\alpha)} \oplus U$.
Let $\aFock = \Fock \ominus \sFock$ denote the antisymmetric part of Fock space.

\begin{lem}\label{L:structure}
Let $U$ be a spherical unitary, and let $L^{(\alpha)} \oplus S$ be the minimal row isometric
dilation of $M^{(\alpha)}\oplus U$ as above. 
Suppose that $Y$ is an operator commuting with $L^{(\alpha)}\oplus S$. 
Write $Y=\begin{sbmatrix}Y_{11} & Y_{12}\\ Y_{21} & Y_{22}\end{sbmatrix}$
with respect to the decomposition $\K' = \Fock^{(\alpha)} \oplus \N'$. Then
\begin{enumerate}
\item $Y_{11}$ commutes with $L^{(\alpha)}_1,\ldots, L^{(\alpha)}_n$.
\item $Y_{22}$ commutes with $S_1,\ldots, S_n$.
\item $Y_{12} = 0$.
\item $P_\N Y_{21}|_{\aFock^{(\alpha)}} = 0$ and 
$U_jP_\N Y_{21}|_{\sFock^{(\alpha)}} = P_\N Y_{21}|_{\sFock)^{(\alpha)}} M_j^{(\alpha)}$
for $1 \le j \le n.$
\end{enumerate}
\end{lem}

\begin{proof}
We have $(L^{(\alpha)}\oplus S)Y^{(n)} = Y(L^{(\alpha)}\oplus S)$, which yields
\[
 \begin{bmatrix}
  L^{(\alpha)}Y_{11}^{(n)} & L^{(\alpha)}Y_{12}^{(n)}\\ SY_{21}^{(n)} & SY_{22}^{(n)} 
 \end{bmatrix} =
 \begin{bmatrix}
  Y_{11}L^{(\alpha)} & Y_{12}S\\ Y_{21}L^{(\alpha)} & Y_{22}S
 \end{bmatrix} .
\]
Thus (i) and (ii) are immediate from the equality of the $11$ and $22$ entries, respectively. 
The $12$ entry yields 
\[ 
 \sum_{i=1}^n L_i^{(\alpha)} Y_{12} S_i^*= L^{(\alpha)}Y_{12}^{(n)} S^* = Y_{12}SS^* = Y_{12} .
\]
Repeated application yields a sum over all words in $\Fn$ of length $k$:
\[ Y_{12} = \sum_{|w|=k} L_w^{(\alpha)} Y_{12} S_w^* .\]
Thus $Y_{12}$ has range contained in the space
$\bigwedge_{k\ge 1} \sum_{|w|=k}(L_{w}\Fock)^{(\alpha)}=\{0\}$. 
Therefore, $Y_{12}=0$; so (iii) holds.

Also, $SY_{21}^{(n)}=Y_{21}L^{(\alpha)}$. 
Let $T=P_\N Y_{21}$. 
Then 
\[
 U_j T = U_jP_\N Y_{21} = P_\N S_jY_{21}
  = P_\N Y_{21}L_j^{(\alpha)} = T L_j^{(\alpha)} .
\] 
Therefore $U_w T = TL_w^{(\alpha)}$ for all words $w \in \Fn$. 
(Here, if $w=i_1\dots i_k$, we set $U_w = U_{i_1} \dots U_{i_k}$.) 
In particular, 
\[
 T L_v(L_iL_j-L_jL_i)L_w = U_v(U_iU_j-U_jU_i)U_wT = 0
\]
for all $1 \le i,j \le n$ and $u,v \in \Fn$.
By \cite[Prop.~2.4]{DP2}, the commutator ideal of $\fL_n$ is spanned by
the words $L_v(L_iL_j-L_jL_i)L_w$, and the range of this ideal is $\aFock$.
Hence $P_\N Y_{21} P_{\aFock^{(\alpha)}} = T P_{\aFock^{(\alpha)}} = 0$;
or $T=TP_{\sFock^{(\alpha)}}$.
Now,
\begin{align*}
 U_{j}T|_{\sFock^{(\alpha)}} & = TL_j^{(\alpha)}|_{\sFock^{(\alpha)}}
  = TP_{\sFock^{(\alpha)}}L_j^{(\alpha)}|_{\sFock^{(\alpha)}} 
  = T|_{\sFock^{(\alpha)}}M_{j}^{(\alpha)}.
\end{align*}
Thus (iv) follows.
\end{proof}

This next lemma is taken from Davie--Jewell \cite[Prop.~2.4]{DavJew}.
We provide a more conceptual proof.

\begin{lem} \label{L:DavJew} 
Let $U = \big[ U_1\ \ldots\ U_n \big]$ be a spherical unitary on a Hilbert space $\H$.
Suppose that $A$ is a bounded operator on $\H$
satisfying $\sum_{j=1}^{n}U_jAU_j^* = A$.
Then $A$ lies in the commutant $\ca(\{U_1,\dots,U_n\})'$.
\end{lem}

\begin{proof}
Let $E$ be the spectral measure for $U$. 
Observe that $E$ is supported on the unit sphere $\bS_n$ of $\bC^n$;
and $U_j = \int_{\bS_n} z_j \,dE$ for $j=1,\ldots,n$. 
We will show that $E(X)AE(Y)=0$ if $X,Y$ are disjoint Borel subsets of $\bS_n$.
Set $B = E(X)AE(Y)$.

Define 
\[ \rho(X,Y) = \inf\{ |1-\langle x,y\rangle| : x\in X, y\in Y \} \]
and  
\[ \diam(X) = \sup\{ |a-b|: a,b\in X \} .\] 
First consider the case in which 
\[  \max \{ \diam(X),\diam(Y) \} < \rho(X,Y)/2n.\] 
Fix points $x=(x_1,\ldots,x_n) \in  X$ and $y=(y_1,\ldots,y_n) \in Y$. 
Since $E(X)$ and $E(Y)$ commute with the $U_j$'s, we have
\begin{align*}
 B &= \sum_{j=1}^{n} U_jBU_j^* = \sum_{j=1}^{n}U_jE(X)BE(Y)U_j^*\\
 &= \sum_{j=1}^{n} x_j\bar{y}_jB + 
 (U_j \!-\! x_j)E(X)BE(Y)U_{j}^* + x_jE(X)BE(Y)(U_{j}^* \!-\!\bar{y}_j) .
\end{align*}
Since $|z_j-x_j| \le \diam(X)$ for all $z\in X$, we have 
\[ \|(U_j-x_j)E(X)\| \le \diam(X) .\] 
Similarly, $\|E(Y)(U_{j}^*-\bar{y}_j)\| \le \diam(Y)$. 
Therefore, 
\begin{align*}
 &|1-\ip{ x,y}| \|B\| \le \\
 &\quad \sum_{j=1}^{n} \|(U_j-x_j) E(X)\|\,\| BE(Y)U_j^* \|
  + \| x_jE(X)B \|\,\|E(Y)(U_j^*-\bar{y}_j) \| \\
 &\quad\leq n\|B\|(\diam(X)+\diam(Y)) < \rho(X,Y) \|B\|.
\end{align*}
But $|1-\langle x,y\rangle|\geq\rho(X,Y)$, 
so $\|B\|=0$.

For $X$ and $Y$ with $\rho(X,Y)>0$, partition $X$ as a disjoint union of Borel subets 
$X_1,\ldots,X_k$ and $Y$ as a disjoint union of Borel subets $Y_1,\ldots,Y_l$ 
so that 
\[
 \max\{\diam(X_i),\diam(Y_j) : 1 \le i \le k,\,1\le j \le l \} < \rho(X,Y)/2n.
\]
Since $ \rho(X_i,Y_j) \ge \rho(X,Y)$, 
it then follows that $E(X_i)AE(Y_j)=0$. 
Thus, $E(X)AE(Y) = \sum_{i,j} E(X_i)AE(Y_j)=0$.
 
Finally suppose $X$ and $Y$ are disjoint Borel subsets of $\bS_n$. 
By the regularity of the spectral measure $E$, we have 
$E(X)=\bigvee E(K)$ as $K$ runs over all compact subsets of $X$. 
For any compact sets $K_1\subset X$ and $K_2\subset Y$, 
we have $\rho(K_1,K_2)>0$ and hence $E(K_1)AE(K_2)=0$. 
It follows that $E(X)AE(Y)=0$. 

In particular $E(X)AE(\bS_n\backslash X) =0 = E(\bS_n\backslash X) A E(X)$.  So
\[ E(X)A = E(X)AE(X) = AE(X) .\] 
Since $A$ commutes with all the spectral projections, we see that $A$ 
commutes with $U_j$ and $U_{j}^*$ for $1\leq j\leq n$, 
and thus lies in $\ca(\{U_1,\dots,U_n\})'$.
\end{proof}

\begin{lem} \label{L:spherical}
Let $U=\big[ U_1\ \ldots\ U_n \big]$ be a spherical unitary on a Hilbert space $\N$. 
Let $S= \big[ S_1\ \ldots\ S_n \big]$ be the minimal isometric dilation of $U$ on $\N'\supset\N$. 
Suppose that $Y$ is an operator on $\N'$ which commutes with $S_1,\ldots, S_n$. 
Then $A=P_\N Y|_\N$ lies in $\ca(\{U_1,\ldots, U_n\})'$ and $P_\N Y P_\N^\perp=0$.
\end{lem}

\begin{proof}
For $1\leq j\leq n$, write $S_j=\begin{sbmatrix}U_j & 0\\ D_j & E_j\end{sbmatrix}$ 
with respect to the decomposition $\N'=\N\oplus\N^\perp$. 
Let us write $D = \big[ D_1\ \ldots\  D_n \big]$ and $E = \big[ E_1\ \ldots\ E_n \big]$. 
Then we may write $S=\begin{sbmatrix}U & 0\\ D & E\end{sbmatrix}$. 
Since $S$ is a row isometry, $S^* S=I_{\N'}^{(n)}$, which implies $U^* U+D^* D=I_\N^{(n)}$ 
(note that $U^*$ and $D^*$ are column operators). 
This together with $UU^* = I_\N$ implies 
\[ UD^* DU^*  = UU^*  - UU^* UU^* =0 .\] 
Thus, $DU^* =0$. 

Since $Y$ commutes with the $S_j$'s, we have $SY^{(n)}=YS$. 
Write $Y = \begin{sbmatrix}A & B\\ * & *\end{sbmatrix}$.
Then $UA^{(n)} = AU+BD$ and $UB^{(n)}=BE$. 
Therefore, 
\[
 \sum_{j=1}^n U_jAU_j^* = UA^{(n)}U^*  = AUU^* + BDU^* = A .
\]
By Lemma \ref{L:DavJew}, $A$ belongs to $\ca(\{U_1,\dots,U_n\})'$.

Since $AU = UA^{(n)} = AU + BD$, we obtain $BD=0$, i.e.\ 
$BD_j=0$ for $1 \le j \le n$. 
The $12$ entry of $SY^{(n)}=YS$ yields $U_j B = B E_j$.
Since $S$ is the minimal dilation, $\N' = \bigvee_{w\in\Fn} S_w \N$. 
Thus $\N^\perp = \bigvee_{w\in\Fn} E_w D \N^{(n)}$.
However $B E_w D \N^{(n)} = U_w BD \N^{(n)} = 0$ for all $w \in \Fn$.
Hence $B=0$.
\end{proof}

We now have the tools to complete the proof of the commutant lifting theorem.

\begin{proof}[\em\textbf{Proof of Theorem~\ref{CLT}}]
We proceed as indicated at the beginning of this section, dilating
$\tilde T \simeq M^{(\alpha)} \oplus U$ to a row isometry on $\K'$
unitarily equivalent to $L^{(\alpha)} \oplus S$.
Let $Y$ be the commuting lifting of $X$, and set $Z = P_\K Y|_\K$.
Since $P_\H Y = X P_\H$, it follows that $P_\H Z = XP_\H$. 
So $Z$ dilates $X$. 
Also, 
\[ \|X\|\leq\|Z\|\leq\|Y\|=\|X\| , \]
and hence $\|Z\|=\|X\|$.
So it remains to verify that $Z$ commutes with $\tilde T_j$ for $1 \le j \le n$.

Write $Y=\begin{sbmatrix} Y_{11} & Y_{12}\\ Y_{21} & Y_{22}\end{sbmatrix}$ 
with respect to the direct sum $\Fock^{(\alpha)}\oplus\N'$. 
By Lemma \ref{L:structure}, $Y_{12}=0$ and $Y_{11}$ commutes with 
$L_{1}^{(\alpha)},\ldots, L_{n}^{(\alpha)}$. 
Since the commutant of $\fL_n$ is $\fR_n$, we see that $Y_{11}$ can be written as an 
$\alpha\times\alpha$ matrix with coefficients in $\fR_n$. 
Since the symmetric Fock space $\sFock$ is co-invariant for $\fR_n$ and $\fL_n$, 
we see that $\sFock^{(\alpha)}$ is co-invariant for $Y_{11}$
and for $L_1^{(\alpha)},\ldots, L_n^{(\alpha)}$. 
Therefore, $P_{\sFock^{(\alpha)}} Y_{11}|_{\aFock^{(\alpha)}}=0$ 
and $P_{\sFock^{(\alpha)}}Y_{11}|_{\sFock^{(\alpha)}}$ commutes with  
$M_j^{(\alpha)}=P_{\sFock^{(\alpha)}}L_{j}^{(\alpha)}|_{\sFock^{(\alpha)}}$ for $1\leq j\leq n$.

By Lemma \ref{L:structure} again, we have 
\[ P_{\N}Y_{21}|_{\aFock^{(\alpha)}} = 0  \]
and
\[
 U_jP_{\N}Y_{21}|_{\sFock^{(\alpha)}}=
 P_{\N}Y_{21}|_{\sFock^{(\alpha)}}M_j^{(\alpha)} 
 \qfor j=1,\ldots,n
\]
and $Y_{22}$ commutes with each $S_j$. 
Lemma \ref{L:spherical} shows that 
$P_{\N}Y_{22}|_{\N^\perp}=0$ and $P_{\N}Y_{22}|_{\N}$ commutes with $U_1,\ldots, U_n$.

Since 
\[
 Z = \begin{bmatrix}
 P_{\sFock^{(\alpha)}}Y_{11}|_{\sFock^{(\alpha)}} & P_{\sFock^{(\alpha)}}Y_{12}|_{\N}\\ 
 P_{\N}Y_{21}|_{\sFock^{(\alpha)}} & P_{\N}Y_{22}|_{\N}
 \end{bmatrix},
\]
$Z$ commutes with 
$M_1^{(\alpha)}\oplus U_1,\ldots, M_n^{(\alpha)}\oplus U_n$ as claimed. 
\end{proof}

We obtain the following corollary by the standard trick of applying the
commutant lifting theorem to $S \oplus T$ and
$\begin{sbmatrix}0&X\\0&0\end{sbmatrix}$.

\begin{cor} \label{C:CLT}
Suppose that $S=\big[ S_1\ \dots\ S_n \big]$ and
$T=\big[ T_1\ \dots\ T_n \big]$ are commuting  row contractions 
on Hilbert spaces $\H_1$ and $\H_2$, respectively, 
and that $X$ is an operator in $\B(\H_2,\H_1)$ such that
$S_i X = X T_i$ for $1 \le i \le n$. 
Let  $\tilde S = \big[ \tilde S_1\ \dots\ \tilde S_n \big]$ and
$\tilde T = \big[ \tilde T_1\ \dots\ \tilde T_n \big]$ be the Arveson dilations
of $S$ and $T$ on $\K_1$ and $\K_2$, respectively. 
Then there is an operator $Z$ in $\B(\K_2,\K_1)$ such that
$\tilde S_i Z = Z \tilde T_i$ for $1 \le i \le n$, 
which dilates $X$ in the sense $P_{\H_1} Z = X P_{\H_2}$
and satisfies $\|Z\|=\|X\|$.
\end{cor}

\section{Failure of General Commutant Lifting}\label{S:fail}

In this section, we consider more general dilations rather than the
canonical co-extension.  Since this is not unique, it should not be
surprising that there is no commutant lifting theorem.

Before doing so, we wish to make the point that this is more 
natural in the single variable case.  The reason is that one
wants to obtain the maximal dilations, i.e. dilations which
can only be further dilated by the adjoining of another direct summand.
These are the elements which determine representations of the underlying
operator algebra which extend to $*$-representations of the enveloping
C*-algebra, and factor through the C*-envelope \cite{DMc}.
In the case of a single contraction, the operator algebra generated 
by a universal contraction is the disk algebra, as evidenced by
(the matrix version of)  the von Neumann inequality.
The C*-envelope of $\AD$ is the abelian algebra $\CT$.
To dilate a contraction to a maximal dilation, one must go to the
unitary dilation.

Note that in the case of a single contraction, the uniqueness of the minimal unitary
dilation means that one can dilate in either direction, or alternate at randon, since
once one obtains a maximal dilation, it will be the minimal unitary dilation
plus a direct summand of some other arbitrary unitary.

However in the case of a commuting row contraction, one is done after
the co-extension to the Arveson dilation.  No further extension is possible
if the maximal co-extension is done first.  So one could argue that the extension, 
being unnecessary and non-canonical, should never be considered.

The following example adds further evidence.

\begin{eg}
Fix $0<r<1$. Let 
\[
 T_1=\begin{bmatrix} r & 0\\ 0 & 0\end{bmatrix} 
 \qand 
 T_2=\begin{bmatrix} 0 & 0\\ 0 & r\end{bmatrix} .
\] 
Clearly $T= \big[ T_1, T_2 \big]$ is a commuting row contraction 
on $\H_1=\bC^2$. 
For $0<\eps\neq r$, define
\begin{equation*} 
 B_1=\begin{bmatrix}\eps & 0 & 0\\ \eps(\eps-r) & r & 0\\ \eps^2 & 0 & 0\end{bmatrix} 
 \qand  
 B_2=\begin{bmatrix}\eps & 0 & 0\\ \eps^2 & 0 & 0\\ \eps(\eps-r) & 0 & r\end{bmatrix}
\end{equation*}
Then $B_1$ and $B_2$ commute and if $\eps$ is small enough, 
$B= \big[ B_1,B_2 \big]$ is a row contraction on $\H_2\oplus\H_1=\bC^3$. 
Let $S= \big[ S_1 \ S_2 \big]$ be the Arveson minimal dilation of $B$, 
acting on $\H_2\oplus\H_1\oplus\H_3$. 
Then $S$ is a maximal dilation of $T$.

Let $X=\begin{sbmatrix}1 & 0\\ 0 & -1\end{sbmatrix}$. 
Observe that $X$ is a unitary operator which commutes with $T_1$ and $T_2$. 
We will show that there does not exist an operator $Y$ on 
$\H_2\oplus\H_1\oplus\H_3$ which is a dilation of $X$ such that 
$\|Y\|=1$ and $Y$ commutes with both $S_1$ and $S_2$. 
Suppose there were such an operator $Y$. 
Write 
\[
 Y = \begin{bmatrix}
  y_{11} & y_{12} & y_{13} & y_{14}\\
  y_{21} & -1 & 0 & y_{24}\\
  y_{31} & 0 & 1  & y_{34}\\
  y_{41} & y_{42} & y_{43} & y_{44}
 \end{bmatrix}
\]
Since $\|Y\|=1$ and $X$ is a unitary, we see that 
\[
 y_{12}=y_{13}=y_{21}=y_{24}=y_{31}=y_{34}=y_{42}=y_{43}=0 .
\]
To be a dilation, $y_{14}=0$ also; but one could also look at the
$24$ entry of $S_1Y-YS_1$ to see that this is necessary.
But then $S_1Y-YS_1$ equals
\begin{gather*}
\begin{bmatrix}
 \eps & 0 & 0 &0\\ 
 \eps(\eps-r) & r & 0 &0\\ 
 \eps^2 & 0 & 0 &0\\
 * & * & * & *
\end{bmatrix}
\begin{bmatrix}
  y_{11} & 0 &0 & 0\\
  0 & -1 & 0 & 0\\
  0 & 0 & 1  & 0\\
  y_{41} & 0 & 0 & y_{44}
 \end{bmatrix}
-
\begin{bmatrix}
  y_{11} & 0 &0 & 0\\
  0 & -1 & 0 & 0\\
  0 & 0 & 1  & 0\\
  y_{41} & 0 & 0 & y_{44}
\end{bmatrix}
\begin{bmatrix}
 \eps & 0 & 0 &0\\ 
 \eps(\eps-r) & r & 0 &0\\ 
 \eps^2 & 0 & 0 &0\\
 * & * & * & *
\end{bmatrix}
\\ =
\begin{bmatrix}
 0 & 0 & 0 & 0\\ 
(y_{11}-1)\eps(\eps-r) & 0 & 0 & 0\\ 
(y_{11}+1)\eps^2 & 0 & 0 & 0\\
 * & * & * & *
\end{bmatrix}
\end{gather*}
It follows that no choice of $y_{11}$ makes both the $21$ and $31$ entries $0$.
Therefore there is no commuting lifting of $X$.
\end{eg}

\section{The Ball--Trent--Vinnikov CLT}\label{S:BTV}

A reproducing kernel Hilbert space (RKHS) is a Hilbert space $\H$ of functions
on a set $X$ which separates points such that point evaluations are continuous.
There are functions $\{k_x \in \H : x \in X\}$ such that $f(x) = \ip{h,k_x}$ 
for all $f\in\H$ and $x\in X$.
The kernel function $k$ on $X \times X$ given by $k(x,y)=k_y(x)$ is positive definite,
meaning that the matrix $\big[k(x_i,x_j)\big]_{i,j=1}^n$ is positive definite for 
every finite subset $x_1,\dots,x_n$ of $X$.
The kernel is called irreducible if $k(x,y) \ne 0$ for all $x,y\in X$.

The space of multipliers $\Mult(\H)$ consists of the (bounded) functions $h$ on $X$
such that $M_h f(x) = h(x) f(x)$ determines a bounded operator.
$\Mult(\H)$ is endowed with the operator norm, and it forms a 
\wot-closed, maximal abelian subalgebra of $\B(\H)$.

The kernel $k$ (or the space $\H$) is called a complete Nevanlinna--Pick kernel 
if for every finite subset $x_1,\dots,x_n$ of $X$ and $d\times d$ matrices
$W_1,\dots,W_n$, there is an operator $F$ in $\fM_d(\Mult(\H))$ such that
$F(x_i) = W_i$ for $1 \le i \le n$ and $\|M_F\| \le 1$ if and only if the Pick matrix
$\big[ (I_d - W_iW_j^*) k(x_i,x_j) \big]_{i,j=1}^n$ is non-negative.
This is motivated by the classical Nevanlinna--Pick Theorem for $H^\infty$,
where the corresponding RKHS is $H^2$ with the Szego kernel 
$k_z(w)=\frac1{1-w\bar{z}}$ for $z \in \bD$.
The necessary and sufficient condition
for the existence of an $\fM_d$-valued $H^\infty$ function $F$ with $F(z_i)=W_i$
for $1 \le i \le n$ and $\|F\|_\infty \le 1$ is precisely that
$\Big[ \frac{I_k - W_iW_j^*}{1-z_i\bar{z_j}}  \Big]$ is non-negative.

These kernels were characterized by McCullough \cite{McC1,McC2} and Quiggin \cite{Q}.
Agler and McCarthy \cite{AM} showed that every irreducible complete Nevanlinna--Pick
kernel is just the restriction of symmetric Fock space on countably many generators
to the span of some subset of the kernel functions.  That is, one can always identify
$X$ with a subset of the ball $\bB_\infty$.

Now we can state the Ball--Trent--Vinnikov theorem \cite{BTV}.

\begin{thm}[Ball--Trent--Vinnikov]
Let $\H$ be an irreducible complete Nevan\-linna--Pick kernel Hilbert space,
and let $\M$ be a subspace of $\H^{(\alpha)}$ which is co-invariant for
$\Mult(\H)^{(\alpha)}$.  
Suppose that $X\in\B(\M)$ commutes with the algebra
$P_\M \Mult(\H)^{(\alpha)}|_\M$.
Then there is an operator $Y$ on $\H^{(\alpha)}$ which commutes with 
$\Mult(\H)^{(\alpha)}$ such that $P_\M Y = X P_\M$ and $\|Y\|=\|X\|$.
\end{thm}

Using the Agler--McCarthy theorem, we identify $\H$ with a coinvariant
subspace of $H^2_n$ spanned by kernel functions, for some $n \le \infty$.
Saying that $X$ commutes with $P_\M \Mult(\H)^{(\alpha)}|_\M$ is equivalent
to saying that it commutes with the commuting row contraction $T$
where $T_i = P_\M M_{z_i}|_\M$ for $1 \le i \le n$.
However, in this case, $T$ is a \textit{pure contraction}, meaning that
the Arveson dilation does not have a spherical unitary part.

In this case, our theorem becomes much easier to prove.
One can dilate $M^{(\alpha)}$ to $L^{(\alpha)}$, and use Popescu's CLT \cite{Pop1} 
to dilate $X$ to an operator $Z$ of the same norm which commutes with
$\fL_n^{(\alpha)}$.  Thus $Z$ belongs to $\fM_\alpha(\fR_n)$.
Restricting $Z$ to the co-invariant subspace $(H^2_n)^{(\alpha)}$ yields
a matrix $Y$ of multipliers of $H^2_n$.  This dilates $X$, has the same norm
and commutes with all multipliers.
Now a further compression to $\H^{(\alpha)}$ yields a matrix of multipliers of $\H$
with these same properties.


\end{document}